%
\newif\ifloadreferences\loadreferencestrue %
%
%
%
%
\let\myfrac=\frac%
\input eplain %
\let\frac=\myfrac%
\input amstex \input epsf %
%
%
\loadeufm\loadmsam\loadmsbm\message{symbol names}\UseAMSsymbols\message{,}%
\magnification 1200 %
\font\myfontdefault=cmr10%
\newif\ifmakebiblio%
\newif\ifinappendices%
\newif\ifundefinedreferences%
\newif\ifchangedreferences%
\makebibliofalse%
\undefinedreferencesfalse%
\changedreferencesfalse%
%
%
%
%
%
\def\setcatcodes{\catcode`\!=0 \catcode`\\=11}%
{\global\let\noe=\noexpand%
\catcode`\@=11 \catcode`\_=11 \setcatcodes%
!global!def!_@@internal@@makeref#1{%
!global!expandafter!def!csname #1ref!endcsname##1{%
!csname _@#1@##1!endcsname%
!expandafter!ifx!csname _@#1@##1!endcsname!relax%
    !write16{#1 ##1 not defined - run saving references}%
    !undefinedreferencestrue%
!fi}}%
!global!def!_@@internal@@makelabel#1{%
!global!expandafter!def!csname #1label!endcsname##1{%
!edef!temptoken{!csname #1info!endcsname}%
!ifloadreferences%
    !expandafter!ifx!csname _@#1@##1!endcsname!relax%
        !write16{#1 ##1 not hitherto defined - rerun saving references}%
        !changedreferencestrue%
    !else%
        !expandafter!ifx!csname _@#1@##1!endcsname!temptoken%
        !else%
            !write16{#1 ##1 reference has changed - rerun saving references}%
            !changedreferencestrue%
        !fi%
    !fi%
!else%
    !expandafter!edef!csname _@#1@##1!endcsname{!temptoken}%
    !edef!textoutput{!write!references{\global\def\_@#1@##1{!temptoken}}}%
    !textoutput%
!fi}}%
!global!def!makecounter#1{!_@@internal@@makelabel{#1}!_@@internal@@makeref{#1}}%
!unsetcatcodes%
}
%
%
%
%
%
\def\turnintolatin#1{\ifcase #1 _\or i\or ii\or iii\or iv\or v\or vi\or vii\or viii\or ix\or x\or xi\or xii\or xiii\or xiv\or xv\or xvi\or xvii\or xviii\or xix\or xx\or xxi\or xxii\or xxiii\or xxiv\or xxv\or xxvi\fi}%
\def\alphanum#1{\ifcase #1 _\or A\or B\or C\or D\or E\or F\or G\or H\or I\or J\or K\or L\or M\or N\or O\or P\or Q\or R\or S\or T\or U\or V\or W\or X\or Y\or Z\fi}%
\newwrite\references%
\ifloadreferences{\catcode`\@=11 \catcode`\_=11 \global\def\_@citation@Aledo{1}
\global\def\_@citation@BallSchGro{2}
\global\def\_@citation@Corro{3}
\global\def\_@citation@Costa{4}
\global\def\_@citation@Galvez{5}
\global\def\_@citation@Huber{6}
\global\def\_@citation@LabA{7}
\global\def\_@citation@LabB{8}
\global\def\_@citation@Osserman{9}
\global\def\_@citation@SmiHPP{10}
\global\def\_@citation@SmiPKS{11}
\global\def\_@proc@ThmMainTheorem{1.2.1}
\global\def\_@subhead@SectionLabouriesCompactnessTheorem{2.1}
\global\def\_@proc@ThmLabouriesCompactnessResult{2.1.1}
\global\def\_@proc@PropShapeOperatorDiverges{2.2.1}
\global\def\_@proc@LemmaContainmentInConvexHull{2.2.2}
\global\def\_@subhead@SubheadHyperbolicCuspEnds{2.3}
\global\def\_@eqn@eqnCuspMetric{\relax \unhbox \voidb@x \hbox {(D)}}
\global\def\_@proc@PropPointedDiskIsContained{2.3.2}
\global\def\_@proc@PropEndPointIsInfinity{2.3.3}
\global\def\_@proc@PropLimitOfImmersion{2.3.4}
\global\def\_@proc@PropMoebiusMapsWhichDiverge{2.4.1}
\global\def\_@proc@PropWeierstrassMapConverges{2.4.2}
\global\def\_@proc@PropSecondPartOfMainResult{2.4.3}
\global\def\_@proc@PropCompleteAndProper{2.5.1}
\global\def\_@proc@PropFirstPartOfMainTheorem{2.5.3}
 }%
\else{\openout\references=references.tex }%
\fi%
%
%
\newcount\headno%
\global\headno=0%
\def\headinfo{\ifinappendices\alphanum\headno\else\the\headno\fi}%
\def\nextheadno{\global\advance\headno by 1 \global\subheadno=0 \global\procno=0 \headinfo}%
\makecounter{head}%
%
%
\newcount\subheadno%
\global\subheadno=0%
\def\subheadinfo{\headinfo.\the\subheadno}%
\def\nextsubheadno{\global\advance\subheadno by 1 \global\procno=0 \subheadinfo}%
\makecounter{subhead}%
%
%
\newcount\procno%
\global\procno=0%
\def\procinfo{\subheadinfo.\the\procno}%
\def\nextprocno{\global\advance\procno by 1 \procinfo}%
\makecounter{proc}%
%
%
\newcount\figno%
\global\figno=0%
\def\figinfo{\subheadinfo.\the\figno}%
\def\nextfigno{\global\advance\figno by 1 \figinfo}%
\makecounter{fig}%
%
%
\newcount\eqnno%
\global\eqnno=0%
\def\eqninfo{\text{(\alphanum{\the\eqnno})}}%
\def\nexteqnno{\global\advance\eqnno by 1 \eqninfo}%
\makecounter{eqn}%
%
%
%
%
%
\def\gobbleeight#1#2#3#4#5#6#7#8{}%
\newcount\citationno%
\global\citationno=0%
\def\citationinfo{\the\citationno}%
\makecounter{citation}%
\newwrite\biblio%
\def\newref#1#2{%
\def\temptext{#2}%
\edef\bibliotextoutput{\expandafter\gobbleeight\meaning\temptext}%
\global\advance\citationno by 1\citationlabel{#1}%
\ifmakebiblio%
    \edef\fileoutput{\write\biblio{\noindent\hbox to 0pt{\hss$[\the\citationno]$}\hskip 0.2em\bibliotextoutput\medskip}}%
    \fileoutput%
\fi}%
\def\cite#1{%
$[\citationref{#1}]$%
\ifmakebiblio%
    \edef\fileoutput{\write\biblio{#1}}%
    \fileoutput%
\fi%
}%
%
%
%
%
\let\mypar=\par%
\edef\Pagetitle={Blank}\headline={\hfil\Pagetitle\hfil}%
\edef\Pagefooter={Blank}\footline={\hfil\Pagefooter\hfil}%
%
%
\newcount\showpagenumflag%
\global\showpagenumflag=0 %
\def\nextoddpage%
{\newpage\ifodd\pageno%
\else\global\showpagenumflag=0 %
\null\vfil\eject%
\global\showpagenumflag=1 %
\fi}%
%
%
\font\headfont=cmb12%
\def\newhead#1%
{\ifhmode\mypar\fi%
\ifnum\headno=0 \else\goodbreak\bigskip\fi%
{\headfont\noindent\nextheadno\ - #1.}
\nobreak\medskip}%
%
%
\def\newsubhead#1%
{\ifhmode\mypar\fi%
\ifnum\subheadno=0 \else\goodbreak\medskip\fi%
{\bf\noindent\nextsubheadno\ - #1.\ }}%
%
%
\newif\ifinproclaim%
\global\inproclaimfalse%
\def\proclaim#1{%
\goodbreak\medskip
\bgroup\inproclaimtrue%
\noindent{\bf #1}%
\nobreak\medskip\sl}%
\def\noskipproclaim#1{%
\goodbreak\medskip%
\bgroup\inproclaimtrue%
\noindent{\bf #1}\nobreak\sl}%
\def\endproclaim{\mypar\egroup\nobreak\medskip\ignorespaces}%
%
%
%
\newcount\xpos\newcount\ypos
\def\makelabelgrid{%
\xpos=-5 \ypos=-5 %
\loop\ifnum\xpos<6 %
{\loop\ifnum\ypos<6 %
\def\labeltext{x}%
\ifnum\xpos=0\def\labeltext{+}\fi%
\ifnum\ypos=0\def\labeltext{+}\fi%
\placelabel[\xpos][\ypos]{\labeltext}%
\advance\ypos by 1 %
\repeat}%
\advance\xpos by 1 %
\repeat}%
\def\placelabel[#1][#2]#3{{%
\setbox10=\hbox{\raise #2cm \hbox{\hskip #1cm #3}}%
\ht10=0pt \dp10=0pt \wd10=0pt \box10}}%
%
%
%
%
\def\myitem#1{\noindent\hbox to .5cm{\hfill#1\hss}}%
%
%
%
%
%
%
%
%
%
\font\sansseriften=cmss10%
\font\sansserifseven=cmss7%
\font\sansseriffive=cmss5%
\newfam\sansseriffam%
\textfont\sansseriffam=\sansseriften%
\scriptfont\sansseriffam=\sansserifseven%
\scriptscriptfont\sansseriffam=\sansseriffive%
\def\mathsf{\fam\sansseriffam}%
%
%
%
\font\boldten=cmb10%
\font\boldseven=cmb7%
\font\boldfive=cmb5%
\newfam\mathboldfam%
\textfont\mathboldfam=\boldten%
\scriptfont\mathboldfam=\boldseven%
\scriptscriptfont\mathboldfam=\boldfive%
\def\mathbf{\fam\mathboldfam}%
%
%
%
\font\mycmmiten=cmmi10%
\font\mycmmiseven=cmmi7%
\font\mycmmifive=cmmi5%
\newfam\mycmmifam%
\textfont\mycmmifam=\mycmmiten%
\scriptfont\mycmmifam=\mycmmiseven%
\scriptscriptfont\mycmmifam=\mycmmifive%
\def\hexa#1{\ifcase #1 0\or 1\or 2\or 3\or 4\or 5\or 6\or 7\or 8\or 9\or A\or B\or C\or D\or E\or F\fi}%
\mathchardef\mathi="7\hexa\mycmmifam7B%
\mathchardef\mathj="7\hexa\mycmmifam7C%
%
%
\font\mymsbmten=msbm10 at 8pt%
\font\mymsbmseven=msbm7 at 5.6pt
\font\mymsbmfive=msbm5 at 4pt%
\newfam\mymsbmfam%
\textfont\mymsbmfam=\mymsbmten%
\scriptfont\mymsbmfam=\mymsbmseven%
\scriptscriptfont\mymsbmfam=\mymsbmfive%
\mathchardef\mybeth="7\hexa\mymsbmfam69%
\mathchardef\mygimmel="7\hexa\mymsbmfam6A%
\mathchardef\mydaleth="7\hexa\mymsbmfam6B%
%
%
%
%
\def\proof{{\noindent\bf Proof:\ }}%
\def\remark{{\noindent\bf Remark:\ }}%
\def\qed{~$\square$}%
\def\makeop#1{\global\expandafter\def\csname op#1\endcsname{{\text{#1}}}}%
\def\makeopsmall#1{\global\expandafter\def\csname op#1\endcsname{{\text{\lowercase{#1}}}}}%
%
%
\def\munion{\mathop{\cup}}%
\def\minter{\mathop{\cap}}%
%
%
\makeop{Ext}%
\makeop{Int}%
\makeop{Dist}%
\makeop{Diam}%
\makeop{Length}%
%
%
%
%
%
\def\mlim{\mathop{{\text{Lim}}}}%
\def\mliminf{\mathop{{\text{LimInf}}}}%
\def\msup{\mathop{{\text{Sup}}}}%
%
%
%
\makeop{Dim}%
\makeop{Ker}%
\makeop{Coker}%
\makeop{Tr}%
\makeop{Adj}%
\makeop{Det}%
\makeop{End}%
\makeop{Lin}%
\makeop{Symm}%
\makeop{Mult}%
%
%
\makeop{dx}%
\makeop{dy}%
\makeop{dz}%
\makeop{dt}%
\makeop{dVol}%
\makeop{dArea}%
\makeop{Supp}%
\makeop{Hess}%
\makeop{Lip}%
%
%
\makeop{Re}%
\makeop{Im}%
\makeop{Arg}%
\makeop{Log}%
\makeop{Exp}%
%
%
\makeopsmall{Cos}%
\makeopsmall{Sin}%
\makeopsmall{Tan}%
\makeopsmall{Sec}%
\makeopsmall{Cosec}%
\makeopsmall{Cot}%
\makeopsmall{ArcCos}%
\makeopsmall{ArcSin}%
\makeopsmall{ArcTan}%
\makeopsmall{ArcSec}%
\makeopsmall{ArcCosec}%
\makeopsmall{ArcCot}%
%
%
\makeopsmall{Cosh}%
\makeopsmall{Sinh}%
\makeopsmall{Tanh}%
\makeopsmall{ArcCosh}%
\makeopsmall{ArcSinh}%
\makeopsmall{ArcTanh}%
%
%
\makeop{Vol}%
\makeop{Area}%
\makeop{Riem}%
\makeop{Ric}%
\makeop{Scal}%
\makeop{Euc}%
\makeop{Imm}%
\makeop{Emb}%
%
%
\makeop{Id}%
\makeop{Ad}%
\makeop{O}%
\makeop{SO}%
\makeop{SL}%
\makeop{GL}%
\makeop{Conf}%
\makeop{Homeo}%
\makeop{Diff}%
\makeop{Isom}%
%
%
\makeop{Ind}%
\makeop{Sig}%
\makeop{Spec}%
%
%
\makeop{Conv}%
\makeop{Max}%
\makeop{Min}%
\makeop{Mod}%
\makeop{Deg}%
\makeop{loc}%
%
%
%
%
%
%
%
%
%
%
%
%
%
 %
%
\newref{Aledo}{Aledo J. A., Espinar J. M., A conformal representation for linear Weingarten surfaces in the de Sitter space, {\sl Journal of Geometry and Physics}, {\bf 57}, no. 8, (2007), 1669--1677}
\newref{BallSchGro}{}
\newref{Corro}{}
\newref{Costa}{}
\newref{Galvez}{G\'alvez J. A., Mart\'\i nez A., Mil\'an F., Complete constant Gaussian curvature surfaces in the Minkowski space and harmonic diffeomorphisms onto the hyperbolic plane, {\sl Tohoku Math. J.}, {\bf 55}, no. 4, (2003), 467--614}
\newref{Huber}{Huber A., On subharmonic functions and differential geometry in the large, {\sl Comment. Math. Helv.}, {\bf 52}, (1957), 13--72}
\newref{LabA}{GAFA}
\newref{LabB}{Lemme de Morse}
\newref{Osserman}{}
\newref{SmiHPP}{}
\newref{SmiPKS}{}
%
%
%
%
\def\Pagetitle{\hfil}
\def\Pagefooter{\hfil{\myfontdefault\folio}\hfil}
\makeop{PSL}%
\makeop{Tanh}%
\makeop{Coth}%
\makeop{Hor}%
\makeop{Ver}%
\makeop{Sys}%
\makeop{Crit}%
\makeop{N}%
\makeop{T}%
\makeop{U}%
\null \vfill
\def\centre{\rightskip=0pt plus 1fil \leftskip=0pt plus 1fil \spaceskip=.3333em \xspaceskip=.5em \parfillskip=0em \parindent=0em}%
\def\textmonth#1{\ifcase#1\or January\or Febuary\or March\or April\or May\or June\or July\or August\or September\or October\or November\or December\fi}
\font\abstracttitlefont=cmr10 at 14pt {\abstracttitlefont\centre
On an Enneper-Weierstrass-type representation of constant Gaussian curvature surfaces in $3$-dimensional hyperbolic space.\par}
\bigskip
{\centre 20th April 2014\par}
\bigskip
{\centre Graham Smith\par}
\bigskip
{\centre Instituto de Matem\'atica,\par
UFRJ, Av. Athos da Silveira Ramos 149,\par
Centro de Tecnologia - Bloco C,\par
Cidade Universit\'aria - Ilha do Fund\~ao,\par
Caixa Postal 68530, 21941-909,\par
Rio de Janeiro, RJ - BRASIL\par}
\bigskip
\noindent{\bf Abstract:\ }For all $k\in]0,1[$, we construct a canonical bijection between the space of ramified coverings of the sphere and the space of complete immersed surfaces in $3$-dimensional hyperbolic space of finite area and of constant extrinsic curvature equal to $k$. We show, furthermore, that this bijection restricts to a homeomorphism over each stratum of the space of ramified coverings of the sphere.
%
%
\par
\vfill
\nextoddpage
\global\pageno=1
\myfontdefault
\def\Pagetitle{\sl On an Enneper-Weierstrass type representation...}
\newhead{Introduction}
\newsubhead{Background}The Enneper-Weierstrass (EW) representation, which describes minimal surfaces in $\Bbb{R}^3$ in terms of one holomorphic function and one holomorphic $1$-form over a given surface, is one of the most remarkable tools used in the study of minimal surfaces, yielding in particular constructions of minimal surfaces satisfying certain unexpected and often surprising properties (c.f. \cite{Costa} and Ch. $8$ of \cite{Osserman}). It is with this motivation that various authors have subsequently studied EW-type representations of other types of surfaces (c.f., for example, \cite{Aledo}, \cite{Corro} and \cite{Galvez}). In this vein, in \cite{SmiHPP}, following the work \cite{LabA} and \cite{LabB} of Labourie, we construct an EW-type representation for immersed surfaces of constant positive extrinsic curvature in $3$-dimensional hyperbolic space, $\Bbb{H}^3$, in terms of one holomorphic function defined over a given surface. More precisely, we showed that for all $k\in]0,1[$, every locally conformal mapping (that is, holomorphic local homeomorphism) from the Poincar\'e disk into the Riemann sphere is the Weierstrass map (defined below) of a unique simply-connected immersed surface in $\Bbb{H}^3$ of constant extrinsic curvature equal to $k$ which is complete in a certain sense. In the current paper, building upon our subsequent results of \cite{SmiPKS} we show that for all $k\in]0,1[$ the Weierstrass map actually defines a homeomorphism between the space of complete immersed surfaces in $\Bbb{H}^3$ of finite area and of constant extrinsic curvature equal to $k$, on the one hand, and the space of ramified coverings of hyperbolic type of the Riemann sphere on the other. We find this result particularly interesting on account of the new research directions it opens up along the frontier between the theory of immersed surfaces on the one hand and Teichm\"uller theory on the other, themes which we propose to investigate in depth in forthcoming work.
\newsubhead{Main result}In order to define the EW-type representation, we first recall the definition of the hyperbolic Gauss map. Let $\opT\Bbb{H}^3$ be the tangent bundle over $\Bbb{H}^3$ and let $\opU\Bbb{H}^3\subseteq\opT\Bbb{H}^3$ be the subbundle of unit vectors over $\Bbb{H}^3$. Let $\partial_\infty\Bbb{H}^3$ be the ideal boundary of $\Bbb{H}^3$ (c.f. \cite{BallSchGro}). The {\bf hyperbolic Gauss map} $\overrightarrow{n}:\opU\Bbb{H}^3\rightarrow\partial_\infty\Bbb{H}^3$ is defined such that for all $X\in\opU\Bbb{H}^3$:
$$
\overrightarrow{n}(X) = \gamma(+\infty) = \mlim_{t\rightarrow+\infty}\gamma(t),\eqnum{\nexteqnno}
$$
\noindent where $\gamma:\Bbb{R}\rightarrow\Bbb{H}^3$ is the unique geodesic such that $\partial_t\gamma(0)=X$. Informally, $\overrightarrow{n}(X)$ is the point in $\partial_\infty\Bbb{H}^3$ towards which $X$ points.
\medskip
\noindent Let $\Sigma:=(i,S)$ be an {\bf immersed surface} in $\Bbb{H}^3$. We recall that this means that $S$ is a surface and $i:S\rightarrow\Bbb{H}^3$ is an immersion. In the sequel, all submanifolds and functions will be taken to be smooth and oriented. Let $N:S\rightarrow\Bbb{H}^3$ be the unit normal vector field over $i$ compatible with the orientation. Let $A\in\Gamma(TS)$ be the shape operator of $i$. We define the {\bf extrinsic curvature} of $i$, $K_i:S\rightarrow\Bbb{R}$ by:
$$
K_i=\opDet(A).\eqnum{\nexteqnno}
$$
\noindent We define the {\bf Weierstrass map} of $i$, $\varphi_i:S\rightarrow\partial_\infty\Bbb{H}^3$ by:
$$
\varphi_i = \overrightarrow{n}\circ N.\eqnum{\nexteqnno}
$$
\noindent We say that $i$ is {\bf locally strictly convex} (LSC) whenever its shape operator is at all points positive definite. Observe that this means that $K$ is everywhere positive. Furthermore, when $i$ is LSC, we chose the orientation such that at every point $N$ points outwards from the convex set bounded locally by the image of $i$ at that point. In this case, it follows from classical hyperbolic geometry that the Weierstrass map is a local homeomorphism (c.f. \cite{BallSchGro}). In particular, since $\partial_\infty\Bbb{H}^3$ carries the conformal structure of the Riemann sphere, upon pulling this conformal structure back through $\varphi_i$, we may suppose that $\varphi_i$ is a locally conformal mapping.
\medskip
\noindent For $k\in]0,1[$, the bijection that we construct may be described schematically as follows:
$$
\left\{\matrix
(i,S)\ \text{s.t.}\hfill\cr
i:S\rightarrow\Bbb{H}^3\ \text{an immersion},\hfill\cr
i\ \text{complete},\hfill\cr
i\ \text{finite area}, \&\hfill\cr
K_i=k.\hfill\cr
\endmatrix\right\}
\longleftrightarrow
\left\{\matrix
(\varphi,\tilde{S}, P)\ \text{s.t.}\hfill\cr
\varphi:S\rightarrow\hat{\Bbb{C}}\ \text{holomorphic},\hfill\cr
\varphi\ \text{non-constant},\hfill\cr
P\subseteq S\ \text{finite},\ \&\hfill\cr
\opCrit(\varphi)\subseteq P.\hfill\cr
\endmatrix\right\}
$$
\noindent Formally, we prove:
\proclaim{Theorem \nextprocno}
\noindent Let $\tilde{S}$ be a compact Riemann surface. Let $P\subseteq\tilde{S}$ be a finite set of points such that $S:=\tilde{S}\setminus P$ is of hyperbolic type. Let $\varphi:\tilde{S}\rightarrow\partial_\infty\Bbb{H}^3$ be a ramified covering with ramification points contained in $P$. Then, for all $k\in]0,1[$, there exists a unique complete LSC immersion $i:S\rightarrow\Bbb{H}^3$ of finite area and of constant extrinsic curvature equal to $k$ such that $\varphi$ is the Weierstrass map of $i$.
\medskip
\noindent Conversely, if $\Sigma:=(i,S)$ is a complete immersed LSC surface in $\Bbb{H}^3$ of finite area and constant extrinsic curvature equal to $k$, for some $k\in]0,1[$, and if $\varphi:S\rightarrow\partial_\infty\Bbb{H}^3$ is the Weierstrass map of $i$, then the Riemann surface $(S,\varphi^*\partial_\infty\Bbb{H}^3)$ is conformally equivalent to a compact Riemann surface $\tilde{S}$ with a finite set $P$ of points removed. Furthermore, $\varphi$ extends to a ramified covering of the sphere by $\tilde{S}$ with ramification points contained in $P$.
\medskip
\noindent Finally, for all $k\in]0,1[$, this mapping defines a homeomorphism between the set of complete immersed LSC surfaces in $\Bbb{H}^3$ of finite area and of constant extrinsic curvature equal to $k$, on the one hand, and the set of pointed ramified coverings of the Riemann sphere, on the other.
\endproclaim
\proclabel{ThmMainTheorem}
\noindent The author is grateful to Asun Jim\'enez Grande, David Dumas, Lucio Rodriguez and Harold Rosenberg for helpful conversations.
\global\advance\headno by 1
\global\subheadno=0
\medskip
\newsubhead{Labourie's compactness theorem} Labourie's compactness theorem (c.f. \cite{LabA}) presents a powerful tool for the study of LSC immersed surfaces of constant extrinsic curvature inside $3$-dimensional manifolds. We first require a few definitions. Let $M:=(M,g)$ be a complete $3$-dimensional Riemannian manifold. Let $TM$ be the tangent space to $M$, let $UM\subseteq TM$ be the unit sphere bundle and let $\pi:UM\rightarrow M$ be the canonical projection. Observe that $\pi$ is distance non-increasing. Let $(i,S)$ be an immersed surface in $M$. Let $N:S\rightarrow UM$ be the unit normal vector field over $\Sigma$ compatible with the orientation. We define the {\bf Gauss lift} $\hat{\Sigma}$ of $\Sigma$ by $(\hat{\mathi},S)$, where $\hat{\mathi}=N$. This terminology merely distinguishes between the mapping $N$, considered as a section of $i^*UM$ over $i$, and the mapping $\hat{\mathi}$, considered as an immersion in its own right into the total space of $UM$. For $k>0$, following \cite{LabA}, we say that $\Sigma$ is a {\bf $k$-surface} whenever $\Sigma$ has constant extrinsic curvature equal to $k$ and $\hat{\Sigma}$ is complete.
\subheadlabel{SectionLabouriesCompactnessTheorem}
\medskip
\noindent Let $\Gamma\subseteq M$ be a complete geodesic. Let $\opN\Gamma\subseteq UM$ be the bundle of unit normal vectors over $\Gamma$. If $(\hat{\mathi},S)$ is an immersed surface in $UM$, then we say that $\hat{\Sigma}$ is a {\bf tube} whenever $\hat{\mathi}$ defines a covering map from $S$ onto $\opN\Gamma$ for some complete geodesic $\Gamma$. Furthermore, we say that $\hat{\Sigma}$ is a tube of {\bf finite order} whenever this covering is of finite order.
\medskip
\noindent Let $(\hat{\mathi}_n,S_n,x_n)$ be a sequence of complete pointed immersed surfaces in $UM$. We say that $(\hat{\mathi}_n,S_n,x_n)$ converges towards the complete pointed immersed surface $(\hat{\mathi}_\infty,S_\infty,x_\infty)$ in the {\bf pointed Cheeger-Gromov} sense whenever there exists a sequence $(\alpha_n)$ of mappings such that:
\medskip
\myitem{(1)} for all $n$, $\alpha_n:S_\infty\rightarrow S_n$;
\medskip
\myitem{(2)} for all $n$, $\alpha_n(x_\infty)=x_n$;
\medskip
\myitem{(3)} for every relatively compact open subset $\Omega\subseteq S_\infty$, there exists $N\in\Bbb{N}$ such that for all $n\geqslant N$, the restriction of $\alpha_n$ to $\Omega$ is a diffeomorphism onto its image; and
\medskip
\myitem{(4)} the sequence $(\hat{\mathi}_n\circ\alpha_n)$ converges in the $C^\infty_\oploc$ sense to $\hat{\mathi}_\infty$.
\medskip
\noindent We refer to $(\alpha_n)$ as a {\bf sequence of convergence maps} for the sequence $(\hat{\mathi}_n,S_n,x_n)$ with respect to the limit $(\hat{\mathi}_\infty,S_\infty,x_\infty)$.
\medskip
\noindent Labourie's compactness theorem is now stated as follows:
\goodbreak
\proclaim{Theorem \nextprocno, {\bf Labourie (1997)}}
\noindent Let $M$ be a complete $3$-dimensional Riemannian manifold. Let $(i_n,S_n,x_n)$ be a sequence of pointed $k$-surfaces in $M$ for some fixed $k>0$. For all $n$, let $\hat{\mathi}_n$ be the Gauss Lift of $i_n$. If there exists a compact set $K\subseteq M$ such that $i_n(x_n)\in K$ for all $n$, then there exists a complete pointed immersed surface $(\hat{\mathi}_\infty,S_\infty,x_\infty)$ towards which $(\hat{\mathi}_n,S_n,x_n)$ subconverges in the pointed Cheeger-Gromov sense. Furthermore, either:
\medskip
\myitem{(1)} $(\hat{\mathi}_\infty,S_\infty)$ is a tube; or
\medskip
\myitem{(2)} $i_\infty:=\pi\circ\hat{\mathi}_\infty$ is an immersion.
\endproclaim
\proclabel{ThmLabouriesCompactnessResult}
\remark Care should be taken in interpreting the limit obtained in the second case of Theorem \procref{ThmLabouriesCompactnessResult}. Indeed, although $\hat{\mathi}_\infty$ is complete, there is no reason to suppose that $i_\infty$ is too. This often leads to counterintuitive phenomena.\qed
\newsubhead{General properties of $k$-surfaces}
\proclaim{Lemma \nextprocno}
\noindent Let $(M,g)$ be a complete $3$-dimensional Riemannian manifold whose isometry group acts co-compactly. Let $(i,S)$ be a $k$-surface in $M$ for some $k>0$. If $(S,i^*g)$ has finite area then for all $B>0$ there exists a compact subset $K\subseteq S$ such that for all $x\in S\setminus K$, $\|A(x)\|>B$.
\endproclaim
\proclabel{PropShapeOperatorDiverges}
\proof Suppose the contrary. There exists $B>0$ and a diverging sequence $(x_n)\in S$ such that $\|A(x_n)\|\leqslant B$ for all $n$. Let $\hat{g}$ be the Sasaki metric over $UM$. For all $n$, let $B_n$ be the unit ball about $x_n$ in $S$ with respect to the metric $\hat{\mathi}^*\hat{g}$. Upon extracting a subsequence, we may suppose that all the $(B_n)$ are disjoint. For all $n$, let $\opArea_n$ be the area of $B_n$ with respect to the metric $i^*g$. We claim that the sequence $(\opArea_n)$ is uniformly bounded below. Indeed, let $K\subseteq M$ be a compact fundamental domain for the isometry group. For all $n$, let $\alpha_n:M\rightarrow M$ be an isometry such that $M_n(i(x_n))\in K$. For all $n$, denote $i_n = \alpha_n\circ i$ and let $\hat{\mathi}_n$ be the Gauss lift of $i_n$. By Labourie's Compactness Theorem, upon extracting a subsequence, there exists a complete pointed immersed surface $(\hat{\mathi}_\infty, S_\infty, x_\infty)$ towards which $(\hat{\mathi}_n,S,x_n)$ converges. Let $(\beta_n)$ be a sequence of convergence mappings for $(\hat{\mathi}_n,S,x_n)$ with respect to the limit $(\hat{\mathi}_\infty,S_\infty,x_\infty)$. Observe that $(i_n\circ\beta_n)$ converges in the $C^\infty_\oploc$ sense to $i_\infty:=\pi\circ\hat{\mathi}_\infty$. However, since $(\|A(x_n)\|)$ is uniformly bounded, $(\hat{\mathi}_\infty,S_\infty)$ is not a tube, and so $i_\infty$ is an immersion. In particular, if $B_\infty$ is the unit ball about $x_\infty$ in $S_\infty$ with respect to the metric $\hat{\mathi}_\infty^*\hat{g}$, and if $\opArea_\infty$ is the area of $B_\infty$ with respect to the metric $i_\infty^*g$, then $\opArea_\infty>0$. Thus:
$$
\mliminf_{n\rightarrow\infty}\opArea_n \geqslant \opArea_\infty > 0.
$$
\noindent In particular, the sequence $(\opArea_n)$ is uniformly bounded below, as asserted. However, since the sequence $(B_n)$ consists of disjoint balls, it follows that $S$ has infinite area. This is absurd, and the result follows.\qed
\medskip
\noindent Now consider the case where $M:=\Bbb{H}^3$ is $3$-dimensional hyperbolic space. For any closed subset $X\subseteq\Bbb{H}^3\munion\partial_\infty\Bbb{H}^3$, let $\opConv(X)$ be its convex hull and let $\overline{B}_r(X)$ be the closure in $\Bbb{H}^3\munion\partial_\infty\Bbb{H}^3$ of the set of all points in $\Bbb{H}^3$ lying at a distance no greater than $r$ from $X$.
\proclaim{Lemma \nextprocno}
\noindent Let $(i,S)$ be a compact LSC immersed surface in $\Bbb{H}^3$ of constant extrinsic curvature equal to $k\in]0,1[$. Then:
$$
i(S)\subseteq\overline{B}_r(\opConv(i(\partial S))),
$$
\noindent where $r=\opTanh^{-1}(k)$.
\endproclaim
\proclabel{LemmaContainmentInConvexHull}
\proof Suppose the contrary. Let $X:=\opConv(i(\partial S))$. Let $x\in S$ be such that $i(x)$ lies at a distance greater than $r$ from $X$. In particular, $x$ is an interior point of $S$. Let $p\in X$ be the closest point to $i(x)$. Let $P$ be the supporting plane to $X$ at $p$ normal to the geodesic joining $p$ to $i(x)$. Let $(P_s)$ be the foliation of $\Bbb{H}^3$ by equidistant planes to $P$ parametrised by (signed) distance from $P$. Let $r':=\msup\left\{s\ |\ i(S)\minter P_s\neq\emptyset\right\}$. Since $S$ is compact, $r'<\infty$ and we may suppose that $i(x)\in P_{r'}$. In particular, $i(S)$ is an interior tangent to $P_{r'}$ at this point. It follows by the geometric maximum principle that the extrinsic curvature of $i$ at $z$ is at least $\opTanh(r')>\opTanh(r)=k$. This is absurd, and the result follows.\qed
\medskip
\newsubhead{Limit points of the immersion} Henceforth, we will assume that $(i,S)$ is a complete finite-area immersed surface in $\Bbb{H}^3$ of constant extrinsic curvature equal to $k$ for some fixed $k\in]0,1[$. We continue to denote by $g$ the metric over $\Bbb{H}^3$, by $\hat{g}$ the Sasaki metric over $\opU\Bbb{H}^3$ and by $\pi:\opU\Bbb{H}^3\rightarrow\Bbb{H}^3$ the canonical projection.
\medskip
\noindent By Gauss' equation, $(S,i^*g)$ has constant {\sl intrinsic} curvature equal to $k-1<0$. Since it has finite area, it follows by H\"uber's Theorem (c.f. \cite{Huber}) that $(S, i^*g)$ is conformally equivalent to a compact surface $\tilde{S}$ with a finite set $P$ of points removed. Furthermore, if we denote by $\Bbb{D}$ the Poincar\'e disk, then for every point $p\in P$ there exists $\epsilon\in]0,1[$ and a conformal mapping $\alpha:\Bbb{D}\rightarrow\tilde{S}$ such that $\alpha(0)=p$ and:
\subheadlabel{SubheadHyperbolicCuspEnds}
$$
(i\circ\alpha)^*g_{ij} = \frac{1}{\left|z\right|^2\opLog(\epsilon\left|z\right|)^2}\delta_{ij},\eqnum{\nexteqnno}
$$
\eqnlabel{eqnCuspMetric}%
\noindent where $\delta$ denotes the Euclidean metric over $\Bbb{D}$. The metric on the right-hand side of \eqnref{eqnCuspMetric} is the standard metric of a finite area hyperbolic cusp.
\medskip
\noindent We henceforth identify $\Bbb{D}^*$ with its image in $S$ and suppress $\alpha$ in what follows. For all $r\in]0,1[$, let $C_r$ be the Euclidean circle of radius $r$ about $0$ in $\Bbb{D}^*$, and let $\opLength(C_r)$ be its length with respect to the cusp metric \eqnref{eqnCuspMetric}. Observe that $(\opLength(C_r))$ tends to $0$ as $r$ tends to $0$.
\proclaim{Proposition \nextprocno}
\noindent There exists a sequence $(r_n)$ converging to $0$ and a point $p_\infty\in\Bbb{H}^3\munion\partial_\infty\Bbb{H}^3$ such that the sequence $(C_{r_n})$ converges to $\left\{p_\infty\right\}$ in the Hausdorff sense.
\endproclaim
\proof By compactness of the family of compact sets, there exists a sequence $(r_n)$ converging to $0$ and a subset $C_\infty$ of $\partial_\infty\Bbb{H}^3\munion\Bbb{H}^3$ such that the sequence $(C_{r_n})$ converges to $C_\infty$ in the Hausdorff sense. Since the sequence $(\opLength(C_{r_n}))$ converges to $0$, $C_\infty$ consists of a single point, and the result now follows.\qed
\medskip
\noindent For all $n$, we henceforth denote $C_n:=C_{r_n}$. For all $n<m$, let $A_{n,m}\subseteq\Bbb{D}^*$ be the anulus bounded by $C_n$ and $C_m$ and for all $n$, let $D_n^*$ be the pointed disk bounded by $C_n$. Let $\hat{\mathi}$ be the Gauss lift of $i$. Observe that since $i^*g$ is complete at $0$, so too is $\hat{\mathi}^*\hat{g}$. For all $n$, denote $X_n:=\overline{B}_r(\opConv(i(C_n)\munion\left\{p_\infty\right\}))$.
\proclaim{Proposition \nextprocno}
\noindent For all $n$, $i(D_n)\subseteq X_n$.
\endproclaim
\proclabel{PropPointedDiskIsContained}
\proof Indeed, for all $n<m$, let $X_{n,m}:=\overline{B}_r(\opConv(i(C_n)\munion i(C_m)))$. By Lemma \procref{LemmaContainmentInConvexHull}, for all $n<m$, $i(A_{n,m})\subseteq X_{n,m}$. For all $n$, $(X_{n,m})$ converges to $X_n:=\overline{B}_r(\opConv(i(C_n)\munion\left\{p_\infty\right\}))$ in the Hausdorff sense as $m$ tends to infinity. The result follows upon taking limits.\qed
\proclaim{Proposition \nextprocno}
\noindent $p_\infty\in\partial_\infty\Bbb{H}^3$.
\endproclaim
\proclabel{PropEndPointIsInfinity}
\proof Suppose the contrary. Let $\hat{\mathi}$ be the Gauss lift of $i$. Let $z_n\in\Bbb{D}^*$ be a sequence converging to $0$. We may suppose that for all $n$, $z_n\in D_n$. Furthermore, since $\hat{\mathi}^*\hat{g}$ is complete at $0$, we may suppose that $z_n$ lies at a distance of at least $n$ from $C_n$ with respect to this metric. In particular, for all $R>0$, there exists $N\in\Bbb{N}$ such that for $n\geqslant N$, $\Omega_{R,n}\subseteq D_n$, where $\Omega_{R,n}$ is the open ball of radius $R$ about $z_n$ in $\Bbb{D}^*$ with respect to the metric $\hat{\mathi}^*\hat{g}$.
\medskip
\noindent For all $n$, by Proposition \procref{PropPointedDiskIsContained}, $i(D_n)\subseteq X_n$. Thus, since $(X_n)$ converges to $\overline{B}_r(p_\infty)$ in the Hausdorff sense as $n$ tends to infinity, for sufficiently large $n$, $i(D_n)\subseteq\overline{B}_{2r}(p_\infty)$. However, by Labourie's Compactness Theorem, there exists a complete pointed immersed surface $(\hat{\mathi}_\infty,S_\infty,z_\infty)$ towards which $(\hat{\mathi},\Bbb{D}^*,z_n)$ subconverges. By Lemma \procref{PropShapeOperatorDiverges}, $(\hat{\mathi}_\infty,S_\infty,z_\infty)$ is a tube. By the preceeding discussion, and taking limits, for all $R>0$, $(\pi\circ\hat{\mathi}_\infty)(\Omega_{R,\infty})\subseteq\overline{B}_{2r}(y_\infty)$, where $\Omega_{R,\infty}$ is the open ball of radius $R$ about $z_\infty$ in $S_\infty$ with respect to the metric $\hat{\mathi}_\infty^*\hat{g}$. In other words, $(\pi\circ\hat{\mathi}_\infty)(S_\infty)\subseteq \overline{B}_{2r}(y_\infty)$. This is absurd, since $(\pi\circ\hat{\mathi}_\infty)(S_\infty)$ is a complete geodesic which is therefore not contained in any compact subset of $\Bbb{H}^3$. The result follows.\qed
\proclaim{Proposition \nextprocno}
\noindent $(i(z))$ tends to $p_\infty$ as $z$ tends to $0$.
\endproclaim
\proclabel{PropLimitOfImmersion}
\proof By Proposition \procref{PropEndPointIsInfinity}, $p_\infty\in\partial_\infty\Bbb{H}^3$. The result now follows by Proposition \procref{PropPointedDiskIsContained} since $(X_n)$ converges to $\left\{p_\infty\right\}$ in the Hausdorff sense as $n$ tends to $\infty$.\qed
\goodbreak
\newsubhead{Limit points of the Weierstrass map}
\proclaim{Proposition \nextprocno}
\noindent Let $(x_n)\in\Bbb{H}^3$ be a sequence converging towards $x_\infty\in\partial_\infty\Bbb{H}^3$. Let $\Gamma\subseteq\Bbb{H}^3$ be a geodesic with end-point $x_\infty$. Let $y_\infty$ be the other end-point of $\Gamma$ and let $z$ be any point of $\Gamma$. If $(\alpha_n)$ is a sequence of isometries of $\Bbb{H}^3$ such that for all $n$, $\alpha_n(x_\infty)=x_\infty$ and $\alpha_n(x_n)=z$, then for every compact subset $K\subseteq\partial_\infty\Bbb{H}^3\setminus\left\{x_\infty\right\}$, $(\alpha_n(K))$ converges to $\left\{y_\infty\right\}$ in the Hausdorff sense as $n$ tends to $\infty$.
\endproclaim
\proclabel{PropMoebiusMapsWhichDiverge}
\proof We identify $\Bbb{H}^3$ with the upper half-space in $\Bbb{R}^3$. Upon applying an isometry, we may suppose that $x_\infty=0$, $y_\infty=\infty$ and $z=(0,0,1)$. For all $n$, let $x_n:=(\xi_n,\eta_n,t_n)$. For all $n$, the mapping $\alpha_n$ is given in these coordinates by:
$$
\alpha_n(x,y,t) = \frac{1}{t_n}(x-\xi_n,y-\eta_n,t).
$$
\noindent The result follows.\qed
\proclaim{Proposition \nextprocno}
\noindent $(\overrightarrow{n}\circ\hat{\mathi})(z)$ tends to $p_\infty$ as $z$ tends to $0$.
\endproclaim
\proclabel{PropWeierstrassMapConverges}
\proof Suppose the contrary. There exists a sequence $(z_n)$ converging to $0$ such that the sequence $(q_n):=((\overrightarrow{n}\circ\hat{\mathi})(z_n))$ converges to $q_\infty\neq p_\infty$. We may suppose that for all $n$, $z_n\in D_n$. Furthermore, since $\hat{\mathi}^*\hat{g}$ is complete at $0$, we may suppose that, for all $n$, $z_n$ lies at a distance of at least $n$ from $C_n$ with respect to this metric. In particular, for all $R>0$, there exists $N\in\Bbb{N}$ such that for all $n\geqslant N$, $\Omega_{R,n}\subseteq D_n$, where $\Omega_{R,n}$ is the open ball of radius $R$ about $z_n$ in $\Bbb{D}^*$ with respect to the metric $\hat{\mathi}^*\hat{g}$. Let $r_\infty$ be any other point of $\hat{\Bbb{C}}$ distinct from both $p_\infty$ and $q_\infty$. Let $\Gamma$ be the geodesic joining $p_\infty$ and $r_\infty$. Let $P$ be a totally geodesic plane in $\Bbb{H}^3$ normal to $\Gamma$. For all $n$, let $\alpha_n$ be an isometry of $\Bbb{H}^3$ such that $\alpha_n(p_\infty)=p_\infty$, $\alpha_n(C_n)\minter\Gamma\neq\emptyset$ and $\alpha_n(i(z_n))\in P$. For all $n$, let $i_n:=\alpha_n\circ i$ and let $\hat{\mathi}_n$ be the Gauss lift of $i_n$.
\medskip
\noindent Let $r:=\opTanh^{-1}(k)$. Since $(\opLength(C_n))$ tends to $0$, we may suppose that, for all $n$, $\opLength(C_n)<r$. Since $\overline{B}_r(\Gamma)$ is convex, it follows that for all $n$, $\opConv(i_n(C_n)\munion\left\{y_\infty\right\})\subseteq \overline{B}_r(\Gamma)$. Thus, by Lemma \procref{LemmaContainmentInConvexHull}, for all $n$, $i_n(D_n)\subseteq\overline{B}_r(\opConv(i_n(C_n))\munion\left\{p_\infty\right\})\subseteq\overline{B}_{2r}(\Gamma)$.
In particular, for all $n$, $i_n(z_n)\in\overline{B}_{2r}(\Gamma)\minter P$, which is a compact set. Thus, by Labourie's Compactness Theorem, there exists a complete pointed immersed surface $(\hat{\mathi}_\infty,S_\infty,z_\infty)$ towards which $(\hat{\mathi}_n,\Bbb{D}^*,z_n)$ converges. Let $(\beta_n)$ be a sequence of convergence maps for $(\hat{\mathi}_n,\Bbb{D}^*,z_n)$ with respect to $(\hat{\mathi}_\infty,S_\infty,z_\infty)$. Observe that $(i_n\circ\beta_n)$ converges towards $i_\infty:=\pi\circ\hat{\mathi}_\infty$ in the $C^\infty_\oploc$ sense. By Lemma \procref{PropShapeOperatorDiverges}, $(\hat{\mathi}_\infty,S_\infty)$ is a tube and so $i_\infty(S_\infty)$ is a complete geodesic. However, for all $R>0$, there exists $N\in\Bbb{N}$ such that for $n\geqslant N$, $i_n(\Omega_{R,n})\subseteq i_n(D_n)\subseteq\overline{B}_{2r}(\Gamma)$. Taking limits, it follows that for all $R>0$, $i_\infty(\Omega_{R,\infty})\subseteq\overline{B}_{2r}(\Gamma)$, where $\Omega_{R,\infty}$ is the open ball of radius $R$ about $z_\infty$ in $S_\infty$ with respect to the metric $\hat{\mathi}_\infty^*\hat{g}$. In other words, $i_\infty(S_\infty)\subseteq\overline{B}_{2r}(\Gamma)$. That is, $i_\infty(S_\infty)$ is a complete geodesic lying at constant distance from $\Gamma$. This geodesic therefore coincides with $\Gamma$. In particular, $\hat{\mathi}_\infty(z_\infty)$ is a unit normal vector to $\Gamma$ at the point $\Gamma\minter P$. That is, $\hat{\mathi}_\infty(z_\infty)$ is tangent to $P$ at this point. Since $P$ is totally geodesic, it follows that $q_\infty:=(\overrightarrow{n}\circ\hat{\mathi}_\infty)(z_\infty)\in\partial_\infty P$. However, $(\alpha_n(q_n))=((\alpha_n\circ\overrightarrow{n}\circ\hat{\mathi})(z_n))=((\overrightarrow{n}\circ\hat{\mathi}_n)(z_n))$ converges to $q_\infty$. Furthermore, by Proposition \procref{PropLimitOfImmersion}, $(i(z_n))$ converges to $p_\infty$ as $n$ tends to $\infty$. Since, by hypothesis, $(q_n)$ remains uniformly bounded away from $p_\infty$, it follows by Proposition \procref{PropMoebiusMapsWhichDiverge} that $(\alpha_n(q_n))$ converges to $r_\infty$ as $n$ tends to $\infty$. In particular, $r_\infty=q_\infty$. This is absurd, and it follows that $(\overrightarrow{n}\circ\hat{\mathi})(z)$ converges to $p_\infty$ as $z$ tends to $0$, as desired.\qed
\medskip
\noindent We now obtain the second part of Theorem \procref{ThmMainTheorem}:
\proclaim{Proposition \nextprocno}
\noindent $(S,\varphi^*\hat{\Bbb{C}})$ is conformally equivalent to a compact Riemann surface $\hat{S}$ with a finite set $P$ of points removed. Furthermore, $\varphi$ extends to a meromorphic map from $\hat{S}$ to $\partial_\infty\Bbb{H}^3=\hat{\Bbb{C}}$. That is, $\varphi$ is a ramified covering.
\endproclaim
\proclabel{PropSecondPartOfMainResult}
\proof Let $\tilde{S}$ and $P$ be as at the beginning of this section. Choose $p\in P$. Let $U$ be a neighbourhood of $p$ homeomorphic to a disk such that $p$ is the only point of $P$ in $U$. By definition $\varphi:=\overrightarrow{n}\circ\hat{\mathi}$. Upon reducing $U$ if necessary, $(U,\varphi^*\hat{\Bbb{C}})$ is conformally equivalent to the annulus $A_c:=\left\{z\ |\ c<\left|z\right|<1\right\}$ where $c\in[0,1[$. Observe that $\varphi$ defines a locally conformal mapping $\tilde{\varphi}:A_c\rightarrow\hat{\Bbb{C}}$. Furthermore, by Proposition \procref{PropWeierstrassMapConverges}, there exists $q\in\hat{\Bbb{C}}$ towards which $\tilde{\varphi}(z)$ converges as $\left|z\right|$ converges to $c$. If $c>0$, then, by classical complex analysis, $\tilde{\varphi}$ is constant. This is absurd, and it follows that $c=0$. Furthermore, by Cauchy's removable singularity theorem, $\tilde{\varphi}$ extends to a holomorphic mapping from $A_0$ into $\hat{\Bbb{C}}$. Since $p\in P$ is arbitrary, the result follows.\qed
\newsubhead{Systoles}The proof of the first part of Theorem \procref{ThmMainTheorem} is a fairly straightforward consequence of the results of \cite{SmiHPP} and \cite{SmiPKS}. Indeed, let $\tilde{S}$ be a compact Riemann surface and let $\varphi:\tilde{S}\rightarrow\hat{\Bbb{C}}$ be a ramified covering of the sphere. Let $P\subseteq\tilde{S}$ be a finite subset containing the ramification points of $\varphi$ such that $\tilde{S}\setminus P$ is of hyperbolic type. Choose $k\in]0,1[$. We identify $\hat{\Bbb{C}}$ with $\partial_\infty\Bbb{H}^3$ in the canonical manner, and, denoting $S:=\tilde{S}\setminus P$, we define $(i_k,S)$ to be the unique $k$-surface in $\Bbb{H}^3$ whose Weierstrass map is $\varphi$. The existence and uniqueness of $\hat{\mathi}_k$ are proven in Theorem $1.4$ of \cite{SmiHPP}.
\goodbreak
\proclaim{Proposition \nextprocno}
\noindent $(S,i_k^*g)$ is complete.
\endproclaim
\proclabel{PropCompleteAndProper}
\proof By Theorem $1.2$ of \cite{SmiPKS}, for all $p\in P$, $i_k(x)$ converges to $\varphi(p)\in\partial_\infty\Bbb{H}^3$ as $x$ converges to $p$. In particular, $i_k$ is proper, and completeness follows.\qed
\medskip
\noindent It remains to show that $(S,i_k^*g)$ has finite area. However, recall that $(S,i_k^*g)$ has constant intrinsic curvature equal to $k-1<0$. It is therefore sufficient to show that for all $p\in P$, $(S,i_k^*g)$ has a hyperbolic cusp at $p$. This is readily proven using systoles. We recall that for all $x\in S$, and for any metric $h$ over $S$, the {\bf systole} of $h$ at $x$, which we denote by $\opSys(h,x)$ is defined to be the length of the shortest topologically non-trivial curve in $S$ passing through $x$.
\proclaim{Proposition \nextprocno}
\noindent There exists $B>0$ such that for all $x\in S$, the systole of $i_k^*g$ at $x$ is no greater than $B$.
\endproclaim
\proof It suffices to prove the result in a neighborhood of every point of $P$. Choose $p\in P$. By Theorem $1.2$ of \cite{SmiPKS}, $\Sigma$ is ``asymptotically tubular'' at $p$. That is, if $\hat{\mathi}_k$ denotes the Gauss Lift of $i_k$, then, near $p$, $\hat{\mathi}_k$ is a graph over a finite order tube (as defined in Section \subheadref{SectionLabouriesCompactnessTheorem}) of a section which converges to zero to every order near infinity. Formally, there exists a tube $(\hat{\mathj},\Bbb{R}\times S^1)$ in $U\Bbb{H}^3$ of finite order; a smooth section $\sigma$ of $\hat{\mathj}^*\opT\opU\Bbb{H}^3$ over $]0,\infty[\times S^1$; a neighbourhood $\Omega$ of $p\in S$ homeomorphic to the disk; and a homeomorphism $\alpha:]0,\infty[\times S^1\rightarrow\Omega\setminus\left\{p\right\}$ such that:
\medskip
\myitem{(1)} $\alpha(s,\theta)$ tends to $p$ as $s$ tends to $+\infty$;
\medskip
\myitem{(2)} for all $(s,\theta)\in ]0,\infty[\times S^1$, $(\hat{\mathi}_k\alpha)(s,\theta) = \opExp(\sigma(s,\theta))$, where $\opExp$ here denotes the exponential map of the total space of $\opU\Bbb{H}^3$; and
\medskip
\myitem{(3)} for all $k\geqslant 0$, $\nabla^k\sigma(s,\theta)$ tends to $0$ as $s\rightarrow +\infty$.
\medskip
\noindent Since $S:=\tilde{S}\setminus P$ is of hyperbolic type, in particular, it is not contractible. Thus, for all $(s,\theta)\in]0,\infty[\times S^1$, $\alpha(\left\{s\right\}\times S^1)$ is a closed curve passing through $\alpha(s,\theta)$ which is topologically non-trivial in $S$. Furthermore, $\opLength(\alpha(\left\{s\right\}\times S^1);\hat{\mathi}_k^*\hat{g})$ converges to $2\pi N$ as $s$ tends to $+\infty$, where $N$ is the order of the tube $(\hat{\mathj},\Bbb{R}\times S^1)$, and $\hat{g}$ denotes the Sasaki metric over $\opU\Bbb{H}^3$. In particular, upon reducing $\Omega$ if necessary, for all $x\in\Omega\setminus\left\{p\right\}$, we have:
$$
\opSys_x(\hat{\mathi}_k^*\hat{g},x) \leqslant 2\pi N+1.
$$
\noindent Since $\pi:\opU\Bbb{H}^3\rightarrow\Bbb{H}^3$ is distance non-increasing, in particular, $i_k^*g=(\pi\circ\hat{\mathi}_k)^*g\leqslant\hat{\mathi}_k^*\hat{g}$. Thus:
$$
\opSys(i_k^*g,x) \leqslant 2\pi N+1.
$$
\noindent Since $p\in P$ is arbitrary, the result follows.\qed
\proclaim{Proposition \nextprocno}
\noindent $(S,i_k^*g)$ has finite area.
\endproclaim
\proclabel{PropFirstPartOfMainTheorem}
\proof Choose $p\in P$. Let $\Omega$ be a neighbourhood of $p$ in $\tilde{S}$ which is homeomorphic to the disk. Since $i_k^*g$ is complete and has constant instrinsic curvature equal to $k-1<0$, $(S,i_k^*g)$ either has a cusp or a funnel at $p$. Since $\opSys(i_k^*g,x)$ is bounded, the singularity at $p$ cannot be a funnel. It is therefore a cusp. In particular, it has finite area. Since $p\in P$ is arbitrary, we conclude that $(S,i_k^*g)$ has finite area as desired.\qed
\medskip
\noindent We now prove Theorem \procref{ThmMainTheorem}:
\medskip
{\bf\noindent Proof of Theorem \procref{ThmMainTheorem}:\ }The first assertion follows from Theorem $1.4$ of \cite{SmiHPP} and Propositions \procref{PropCompleteAndProper} and \procref{PropFirstPartOfMainTheorem}. The second assertion follows from Proposition \procref{PropSecondPartOfMainResult}. The final assertion follows from Theorem $1.5$ of \cite{SmiHPP}, and this completes the proof.\qed
\goodbreak
\newhead{Bibliography}
{\leftskip = 5ex \parindent = -5ex
\leavevmode\hbox to 4ex{\hfil \cite{Aledo}}\hskip 1ex{Aledo J. A., Espinar J. M., A conformal representation for linear Weingarten surfaces in the de Sitter space, {\sl Journal of Geometry and Physics}, {\bf 57}, no. 8, (2007), 1669--1677}
\medskip
\leavevmode\hbox to 4ex{\hfil \cite{BallSchGro}}\hskip 1ex{Ballmann W., Gromov M., Schroeder V., {\sl Manifolds of nonpositive curvature}, Progress in Mathematics, {\bf 61}, Birkh\"auser Boston Inc., Boston, MA, (1985)}
\medskip
\leavevmode\hbox to 4ex{\hfil \cite{Corro}}\hskip 1ex{Corro A. V., Mart\'\i nez A., Mil\'an F., Complete flat surfaces with two isolated singularities in hyperbolic $3$-space, {\sl J. Math. Anal. Appl.}, {\bf 366}, (2010), no. 2, 582--592}
\medskip
\leavevmode\hbox to 4ex{\hfil \cite{Costa}}\hskip 1ex{Costa C. J., Example of a complete minimal immersion in $\Bbb{R}^3$ of genus one and three embedded ends, {\sl Bol. Soc. Bras. Mat.}, {\bf 15}, nos. 1-2, (1984), 47--54}
\medskip
\leavevmode\hbox to 4ex{\hfil \cite{Galvez}}\hskip 1ex{G\'alvez J. A., Mart\'\i nez A., Mil\'an F., Complete constant Gaussian curvature surfaces in the Minkowski space and harmonic diffeomorphisms onto the hyperbolic plane, {\sl Tohoku Math. J.}, {\bf 55}, no. 4, (2003), 467--614}
\medskip
\leavevmode\hbox to 4ex{\hfil \cite{Huber}}\hskip 1ex{Huber A., On subharmonic functions and differential geometry in the large, {\sl Comment. Math. Helv.}, {\bf 52}, (1957), 13--72}
\medskip
\leavevmode\hbox to 4ex{\hfil \cite{LabA}}\hskip 1ex{Labourie F., Probl\`emes de Monge-Amp\`ere, courbes pseudo-holomorphes et laminations (French), {\sl G.A.F.A.}, {\bf 7}, (1997), 496--534}
\medskip
\leavevmode\hbox to 4ex{\hfil \cite{LabB}}\hskip 1ex{Labourie F., Un lemma de Morse pour les surfaces convexes (French), {\sl Invent. Math.}, {\bf 141}, no. 2, (2000), 239--297}
\medskip
\leavevmode\hbox to 4ex{\hfil \cite{Osserman}}\hskip 1ex{Osserman R., {\sl A survey of minimal surfaces}, Courier Dover Publications, (2002)}
\medskip
\leavevmode\hbox to 4ex{\hfil \cite{SmiHPP}}\hskip 1ex{Smith G., Hyperbolic Plateau problems, to appear in {\sl Geom. Dedicata}}
\medskip
\leavevmode\hbox to 4ex{\hfil \cite{SmiPKS}}\hskip 1ex{Smith G., Pointed $k$-surfaces, {\sl Bull. Soc. Math. France}, {\bf 134}, no. 4, (2006), 509--557}
\par}
\enddocument
%
%
%
%
\enddocument